\documentclass[12pt]{amsart}
\usepackage{etex,url}
\makeindex

\usepackage{graphicx,tikz,url,hyperref}
\usepackage[latin1]{inputenc} \usepackage[T1]{fontenc} \usepackage{lmodern}
\usepackage{pgf}
\usepackage{epsfig,graphicx,color,amsmath,a4wide,amssymb,amsfonts,amsbsy,xy,latexsym,epsf,pstricks,verbatim,times}
\usepackage{MnSymbol}
\usepackage{pgf,color}

\usepackage{changebar}

\definecolor{LightGrey}{rgb}{.85,.85,.85}
\definecolor{DarkGrey}{rgb}{.5,.5,.5}

\definecolor{Blue}{rgb}{.0,.0,0.9}
\definecolor{LightBlue1}{rgb}{.2,.4,0.9}
\definecolor{LightBlue2}{rgb}{.3,.5,0.9}
\definecolor{LightBlue3}{rgb}{.4,.6,0.9}
\definecolor{LightBlue4}{rgb}{.5,.7,.9}
\definecolor{LightBlue5}{rgb}{.6,.8,.9}
\definecolor{LightBlue6}{rgb}{.7,.9,.9}

\definecolor{Red}{rgb}{.9,.0,.0}
\definecolor{LightRed1}{rgb}{0.9,.2,.4}
\definecolor{LightRed2}{rgb}{0.9,.3,.5}
\definecolor{LightRed3}{rgb}{0.9,.4,.6}
\definecolor{LightRed4}{rgb}{.9,.5,.7}
\definecolor{LightRed5}{rgb}{.9,.6,.8}
\definecolor{LightRed6}{rgb}{.9,.7,.9}

\def\bK{{\mathbf K}}
\def\bM{{\mathbf M}}

\def\bZ{{\mathbf Z}}
\def\bKb{{\bar {\mathbf K}}}
\def\bL{{\mathbf L}}
\def\bZ{{\mathbf Z}}

\def\cO{\mathcal O}
\def\cB{\mathcal B}

\xyoption{all}
\usepackage[english]{babel}
\newcounter{noalgo}[section]
\newdimen\indentalgo
\newdimen\indentalgodec\indentalgo=0.0mm\indentalgodec=10mm

\def\<<{\leavevmode
  \raise0.28ex\hbox{$\scriptscriptstyle\langle\!\langle$}\nobreak
  \hskip -.6pt plus.3pt minus.2pt\,}
\def\>>{\,\nobreak\hskip -.6pt plus.3pt minus.2pt
  \raise0.28ex\hbox{$\scriptscriptstyle\rangle\!\rangle$}}

\def\rgot{{\mathfrak r}}

\newcommand{\mmu}{{\boldsymbol{\mu}}}
\def\rank{\mathop{\rm{rank}}\nolimits }

\def\Hom{\mathop{\rm{Hom}}\nolimits }
\def\Bil{\mathop{\rm{Bil}}\nolimits }
\def\Ker{\mathop{\rm{Ker}}\nolimits }

\def\cQ{{\mathcal Q}}

\newtheorem{theorem}{Theorem}

\providecommand{\myproofname}{Proof}

\usepackage{enumitem}

\usepackage{scalerel,stackengine}
\stackMath
\newcommand\widecheck[1]{%
\savestack{\tmpbox}{\stretchto{%
  \scaleto{%
    \scalerel*[\widthof{\ensuremath{#1}}]{\kern-.6pt\bigwedge\kern-.6pt}%
    {\rule[-\textheight/2]{1ex}{\textheight}}
  }{\textheight}%
}{0.5ex}}%
\stackon[1pt]{#1}{\scalebox{-1}{\tmpbox}}%
}

\begin{document}

\begin{abstract}
  We study the complexity of multiplication of two elements 
  in  a finite field extension
  given by their coordinates in a normal basis.
  We show how to control
  this complexity using the arithmetic
  and geometry of  algebraic curves.
  
\end{abstract}

\title{The equivariant complexity of multiplication in finite
  field extensions}

\author{Jean-Marc Couveignes}
\address{Jean-Marc Couveignes, Univ. Bordeaux, CNRS, INRIA, 
  Bordeaux-INP, IMB, UMR 5251, F-33400 Talence, France.}
\email{Jean-Marc.Couveignes@u-bordeaux.fr}
\author{Tony Ezome}
\address{Tony Ezome, Universit{\'e} des Sciences et Techniques de Masuku,
Facult{\'e} des Sciences, D{\'e}partement de math{\'e}matiques et informatique,
BP 943 Franceville, Gabon.}
\email{tony.ezome@gmail.com}

\date{\today}

\maketitle
\setcounter{tocdepth}{1} 
\tableofcontents

\section{Introduction}

Let $\bK$ be a finite field of  cardinality $q$.
Let $\bL/\bK$ be a finite field extension
of degree $n$.
Given a normal $\bK$-basis $\cB$ of $\bL$
we can represent
elements in $\bL$ by their coordinates
in $\cB$. Exponentiation by $q$ then  corresponds
to a cyclic shift of coordinates and can be computed at almost
no cost. It is a natural concern in this context  to bound the computational complexity of computing
the product of two elements of $\bL$ given by their
coordinates in $\cB$. There is a rich litterature about constructing
normal bases where the cost of multiplication is as small as possible.
See \cite{GT} for  a survey.
In this work we 
 define and study
the  symmetric equivariant
complexity $\nu_q^{sym}(n)$ of multiplication in the  finite field extension $\bL/\bK$.
This is the  Galois equivariant counterpart to the symmetric bilinear 
complexity $\mu_q^{sym}(n)$. It is the size of the smallest decomposition of the multiplication
tensor as a sum of pure  equivariant tensors.  This is an invariant
of the field extension $\bL/\bK$ in the sense that it only
depends on $q$ and $n$. 
While
the  symmetric bilinear 
complexity $\mu_q^{sym}(n)$ partially  controls the cost of multiplication in $\bL$ (it only accounts for bilinear operations), in contrast,  the  symmetric equivariant
complexity $\nu_q^{sym}(n)$ provides an  asymptotic  estimate for the total cost of multiplication
in any normal basis:
the linear part of the calculation consists of $3\nu_q^{sym}(n)$ convolution products, each of them
beeing  computed at the expense of
$O(n\log (n) |\log (\log (n))|)$ operations in $\bK$.
We are interested in proving upper bounds for 
$\nu_q^{sym}(n)$. For example we  prove that $\nu_q^{sym}(n)$ is bounded by
a constant times $\lceil \log_q n\rceil$ in full generality.
This implies that
multiplication in any normal basis
requires
no more than $n(\log n)^{2+o(1)}$ operations in $\bK$.
We also provide methods to bound
 $\nu_q^{sym}(n)$ for given $q$ and $n$.

Section~\ref{sec:tour} is a quick tour of various definitions of complexity in the context
of multiplication in finite field extensions.
In Section~\ref{sec:ac} we recall the elementary properties of the algebraic
complexity of a bilinear map. We introduce in Section~\ref{sec:eqac} the equivariant complexity
of a $C$-equivariant bilinear map, where $C$ is a given finite group.
We prove in Section~\ref{sec:ub}
an inequality between the  equivariant complexity of a $C$-equivariant bilinear map
and the bilinear  complexity of its coordinates.
Sections~\ref{sec:cff} and~\ref{sec:trunc} recall classical results
about the bilinear complexity of multiplication in finite field extensions and truncated power series algebras.
The Galois equivariant complexity of multiplication
in  finite  field extensions is  introduced in Section~\ref{sec:ecff}.
A useful generalization to semi-simple
algebras is introduced in Section~\ref{sec:semi}. The effect of extension
and restriction of scalars on (equivariant) complexities is studied in Sections~\ref{sec:ext1} and \ref{sec:ext2}.
We present in Section~\ref{sec:geo}, \ref{sec:heu}, and
\ref{sec:NS} a general geometric recipe to bound from above the Galois equivariant
complexity of multiplication in a finite  extension $\bL/\bK$ of finite fields.
We first 
construct a cyclic  cover
$\rho : Y\rightarrow X$ between two $\bK$-curves,  then 
realise $\bL/\bK$ as the   residual algebra
of the fiber of $\rho$ above some rational
point on $X$. Evaluation and interpolation on $Y$ naturally produce $\bK[C]$-linear maps.
In the special case when
$X$ and $Y$ are elliptic curves
our construction
generalizes the one presented in \cite{cel}.
The  Chudnovsky's method \cite{CC,STV,BR,CHA,RAN} to bound $\mu_q^{sym}(n)$
relies on the existence of
families
of curves  having an increasing number of rational points while
the  genus is bounded by a constant times the number
of points.
Our construction requires
Jacobians of smallest possible
dimension having a point of given order.
In Sections~\ref{sec:ec} and \ref{sec:gb} we enhance   the specific case when $Y$ and $X$ both have genus one.
Although this special case is not optimal (because we lack  rational points
on elliptic curves when $q$ is small compared to $n$) we have enough control
on the group of points on an  elliptic curve to  prove a satisfactory asympotic statement,
using the general properties of equivariant complexity established in Sections \ref{sec:ub}, \ref{sec:ecff}, \ref{sec:ext2}.
In  Section~\ref{sec:exp} we explain  
how  to better bound  $\nu_q^{sym}(n)$ for  given $q$  and $n$ using the
construction of Section \ref{sec:geo}. 
We experiment   with three examples in Sections
\ref{sec:5}, \ref{sec:239}, \ref{sec:4639}. These examples
illustrate  how the
knowledge of special linear series on low genus curves helps
bounding  $\nu_q^{sym}(n)$ at a minimal   computational cost.
We conclude in Section \ref{sec:concl} with  remarks and questions.

This study has been carried out with financial support from 
the French State, managed by
CNRS in the frame
of the {\it Dispositif de Soutien aux Collaborations avec l'Afrique subsaharienne}
and by 
the French National Research Agency (ANR) 
in the frame of the  Programmes
CIAO (ANR-19-CE48-0008), 
FLAIR  (ANR-17-CE40-0012
ANR-10-IDEX-03-02) and CLapCLap (ANR-18-CE40-0026), and by the
Simons foundation. 

Experiments presented in this paper were carried out using the
PlaFRIM experimental testbed, supported by Inria, CNRS (LABRI and IMB),
Universit{\'e} de Bordeaux, Bordeaux INP and Conseil R{\'e}gional d'Aquitaine
(see \url{https://www.plafrim.fr/}).

We thank the referee for his useful comments.

\section{Various complexities}\label{sec:tour}

There are several notions of complexity in the context
of multiplication in a degree $n$ extension $\bL/\bK$ of  finite  fields.
We assume that  we are given a basis $\cB$ and the coordinates
of the two operands in $\cB$. The output consists of the coordinates of
the product in the basis $\cB$.

In the computational model of
straight line programs, one may
count all arithmetic operations in $\bK$ :
additions, subtractions, multiplications. 
Another option is to  omit additions,
subtractions, and multiplications by a constant in $\bK$. One then
only counts multiplications  of two registers.
This can be justified  if the number of additions, subtractions
and multiplications by a constant, is of the same order of
magnitude as the number of multiplications.

In a more algebraic setting one may count the non-zero
coordinates  of the multiplication tensor in the
basis $\hat\cB \otimes \hat \cB \otimes \cB$. When
$\cB$ is a normal basis, this number can be written
$n\times C_\cB$ where $C_\cB$ is an integer often called
the complexity of the normal basis $\cB$.
It was shown by Mullin, Onyszchuk, Vanstone and
Wilson  \cite{MOVW} that $C_\cB$
is at least $2n-1$. This means that if we only allow products between
the coordinates of the inputs (no intermediate result) the number
of arithmetic operations is at least quadratic in $n$, and most of the time
even cubic. This is a rather pessimistic model that
is well adapted to low capacity computing devices.

A more intrinsic algebraic approach is to define the bilinear complexity
of multiplication in $\bL/\bK$ as the rank $r$ of the multiplication
tensor. The rank is independent of the basis. Given a decomposition
of the multiplication tensor as a sum of $r$ pure tensors, we can compute
products at the expense of $r$ multiplications
between two registers,  $3rn$ multiplications by a constant, and
$3r(n-1)$ additions. According to Chudnovsky and Chudnovsky,
 $r$ is bounded by a constant times $n$. But this says little 
about the  cost of the linear part of the algorithm, since the bound
$3rn$ is  quadratic in $n$. 

We define in Sections \ref{sec:eqac} and \ref{sec:ecff} the equivariant algebraic complexity
of multiplication in $\bL/\bK$. The underlying idea is to stick to the intrinsic algebraic
approach but restrict
the linear part of the algorithm
to Galois equivariant linear forms : convolution products in the algebra of the Galois group.
Respecting the symmetries of the problem is a natural
restriction
in view of the
importance of convolution
products
in fast arithmetic
. See \cite{heid, gathen,GGPold}.

\section{Algebraic complexity of a bilinear map}\label{sec:ac}

We  
recall standard definitions about  complexity
of bilinear maps.
A complete introduction can be found in \cite{BCS}[Chapter 14].
Let $\bK$ be a commutative field. Let $V$ and $W$ be two finite
dimensional $\bK$-vector spaces. Let
\[t : V\times V\rightarrow W\] be a $\bK$-bilinear map.
 We let $\hat V$ be the dual of $V$.
 For $\phi_1$, $\phi_2$ in $\hat V$ and $w$ in $W$ we define the bilinear
 map
  \begin{equation*}
\xymatrix@R-2pc{
  \pi_{w,\phi_1,\phi_2}  & \relax : &    V\times V    \ar@{->}[r]  & W\\
  &&(v_1,v_2)   \ar@{|->}[r] & \phi_1 (v_1)\phi_2 (v_2)w
}
  \end{equation*}
  \smallskip
and we say that $\pi_{w,\phi_1,\phi_2}$ is a {\bf pure}  bilinear map.

If $\phi_1=\phi_2=\phi$ we write  $\pi_{w,\phi}$ for
 $\pi_{w,\phi,\phi}$ and call $\pi_{w,\phi}$ a {\bf pure symmetric}  bilinear map.
For $t$ a 
 $\bK$-bilinear map we define the  {\bf bilinear  complexity}
$R_\bK(t)$  of $t$ to be the smallest integer  such   that
$t$ is the sum of $R_\bK(t)$ pure  bilinear maps.
 In case  $t$ is symmetric we define the {\bf symmetric 
 complexity} $S_\bK(t)$ of $t$ to be the smallest integer  such   that
 $t$ is the sum of $S_\bK(t)$ pure symmetric bilinear
 maps.
  Equivalently $S_\bK(t)$
 is the smallest integer $k$ such that there exist
 two $\bK$-linear maps
 \[\top : V\rightarrow \bK^{k} \text{ \,\, and \,\,}
 \bot : \bK^{k} \rightarrow  W\]
 such that \[t(l_1,l_2)=\bot (\top (l_1) \bullet_k  \top (l_2))\]
 where the $\bullet_k$ between $\top (l_1)$ and $\top (l_2)$ stands
 for the componentwise product in $\bK^k$.

The vector space of bilinear maps  has
a basis consisting of pure bilinear maps. So any bilinear map  $t : V\times V \rightarrow W$
has complexity \[R_\bK(t)\leqslant \dim W \times (\dim V)^2.\]
The vector space of symmetric bilinear maps has
a basis consisting of $\dim V \times \dim W$ pure symmetric
bilinear maps and $\dim V \times (\dim V -1)/2\times \dim W$ maps of the form
\[\pi_{w,\phi_1,\phi_2}+\pi_{w,\phi_2,\phi_1} = \pi_{w,\phi_1+\phi_2}-\pi_{w,\phi_1}-
\pi_{w,\phi_2}.\]
So any symmetric bilinear map $t : V\times V \rightarrow W$
has symmetric complexity \[S_\bK(t)\leqslant \dim W \times (\dim V)\times (3\dim V-1)/2.\]


 If $t_1 : V_1\times V_1\rightarrow W_1$  and $t_2 : V_2\times V_2\rightarrow W_2$ are
 two symmetric $\bK$-bilinear map, we say that $t_2$ is a {\bf restriction}
 of $t_1$ if there exist
 two $\bK$-linear maps $\top : V_2\rightarrow V_1$
 and $\bot : W_1\rightarrow W_2$ such that
 $t_2=\bot \circ t_1 \circ (\top \times \top)$. It follows that
 $S_\bK(t_2)\leqslant S_\bK(t_1)$. In case the maps
 $\top$ and $\bot$ are bijective we say that $t_1$ and $t_2$
 are {\bf isomorphic}.

\section{Equivariant algebraic complexity}\label{sec:eqac}

Let $C$ be a finite  group of order $n$.
Let $\bK$ be a commutative field. Let $\bK [C]$ be the group algebra. We
denote 
 \begin{equation*}
\xymatrix@R-2pc{
  \star  & \relax : & \bK [C]\times \bK [C]     \ar@{->}[r]  & \bK [C]\\
  &&(\displaystyle \sum_{c\in C}a_c.c,  \displaystyle \sum_{c\in C}b_c.c)\ar@{|->}[r]  & \displaystyle \sum_{c\in C} \,\,
\displaystyle  \sum_{\stackrel{c_1, c_2 \, \in C}{c_1c_2=c}}(a_{c_1}b_{c_2}).c
}
\end{equation*}
 the (convolution) product in $\bK [C]$.
 Considering the coefficients
 $(a_c)_{c\in C}$ in
$\sum_{c\in C}a_c.c$
 as a map
$a : C\rightarrow \bK$ we obtain a natural
isomorphism  of $\bK$-vector spaces
 \begin{equation*}
\xymatrix@R-2pc{
   \bK [C]     \ar@{->}[r]  &\Hom (C, \bK) \\
   \sum_{c\in C}a_c.c \ar@{|->}[r] &  \, (  c\mapsto a_c ).
}
\end{equation*}
 between the group algebra $\bK[C]$ and
  the algebra  $\Hom (C, \bK)$ of maps from $C$ to
 $\bK$.
 Through this identification
 the group algebra inherits a componentwise product
 \begin{equation*}
\xymatrix@R-2pc{
  \diamond  & \relax : & \bK [C]\times \bK [C]     \ar@{->}[r]  & \bK [C]\\
  &&(\sum_{c\in C}a_c.c, \sum_{c\in C}b_c.c)\ar@{|->}[r]  & \sum_{c\in C}(a_cb_c).c
}
\end{equation*}

 For any positive integer
 $k$ we denote $\diamond_k$ the map
\begin{equation*}
\xymatrix@R-2pc{
  \diamond_k  & \relax : & \left(\bK [C]\right)^k \times \left(\bK [C]\right)^k     \ar@{->}[r]  & \left(\bK [C]\right)^k\\
  &&((a_i)_{1\leqslant i\leqslant k}, (b_i)_{1\leqslant i\leqslant k} )\ar@{|->}[r]  &
  (a_i\diamond b_i)_{1\leqslant i\leqslant k}
}
\end{equation*}

\bigskip

  If  $L$ and $M$ are two
 finitely generated left $\bK[C]$-modules,
we say that a  $\bK$-bilinear
 map \[t  : L\times L \rightarrow M\]  is
 a  $C$-{\bf equivariant bilinear map}
 if \[t(c \cdot l_1, c \cdot l_2)=c\cdot t(l_1,l_2)\] for any $l_1$, $l_2$ in $L$
 and $c$ in $C$.
 If
 $\alpha_1$ and  $\alpha_2$ are two $\bK [C]$-linear maps from
 $L$ to $\bK [C]$, and if 
 $m$ is a vector in $M$, we define the $C$-equivariant $\bK$-bilinear
 map
 \begin{equation*}
 \xymatrix@R-2pc{
     \gamma_{m,\alpha_1,\alpha_2} & \relax : &  L\times L \ar@{->}[r] & M \\
     & & (l_1,l_2)  \ar@{|->}[r]& (\alpha_1(l_1)\diamond \alpha_2(l_2)).m
   }
 \end{equation*}
 We say that $\gamma_{m,\alpha_1,\alpha_2}$ is
 a {\bf pure} $C$-{\bf equivariant} $\bK$-bilinear map.
If $\alpha_1=\alpha_2=\alpha$ we write  $\gamma_{m,\alpha}$ for
$\gamma_{m,\alpha,\alpha}$
and call $\gamma_{m,\alpha}$ a {\bf pure symmetric}  $C$-{\bf equivariant}
$\bK$-bilinear map.
For $t$ a $C$-equivariant
 $\bK$-bilinear map we define the {\bf equivariant complexity}
 of $t$ to be the smallest integer $R_{\bK,C}(t)$ such   that
 $t$ is the sum of $R_{\bK,C}(t)$ pure $C$-equivariant maps.
 In case  $t$ is symmetric we define the {\bf symmetric equivariant
 complexity} of $t$ to be the smallest integer $S_{\bK,C}(t)$ such   that
 $t$ is the sum of $S_{\bK,C}(t)$ pure symmetric
 $C$-equivariant $\bK$-bilinear maps.
 Equivalently $S_{\bK,C}(t)$
 is the smallest integer $k$ such that there exist
 two $\bK[C]$-linear maps
 \[\top  : L\rightarrow \left(\bK[C]\right)^{k} \text{ \,\, and \,\,}
 \bot  : \left(\bK[C]\right)^{k} \rightarrow  M\]
 such that \[t(l_1,l_2)=\bot (\top (l_1)\diamond_k \top (l_2)).\]

 \bigskip
 

 If $t_1 : L_1\times L_1\rightarrow M_1$  and $t_2 : L_2\times L_2\rightarrow
 M_2$ are
 two symmetric $\bK$-bilinear $C$-equivariant maps, we say that $t_2$ is a {\bf restriction}
 of $t_1$ if there exist
 two $\bK[C]$-linear maps $\top : L_2\rightarrow L_1$
 and $\bot : M_1\rightarrow M_2$ such that
 $t_2=\bot \circ t_1 \circ (\top \times \top)$. It follows that
 $S_{\bK, C}(t_2)\leqslant S_{\bK, C}(t_1)$. In case the maps
 $\top$ and $\bot$ are bijective we say that $t_1$ and $t_2$
 are {\bf isomorphic}.

\section{General upper bounds}\label{sec:ub}

Let $C$ be a finite group
of order $n$.
Let $e$ be the identity element in
 $C$. Let $M$ be a left $\bK[C]$-module.
 We let \[\hat M = \Hom_{\bK}(M,\bK)\] be the dual of $M$
 as a $\bK$-vector space. Let
 \[\check M =  \Hom_{\bK[C]}(M,\bK[C])\] be the dual of $M$ as a $\bK[C]$-module.
 For any $\phi$ in  $\check M$ and $m$ in $M$ we write
 \[\phi (m)= \sum_{c\in C} \phi_c (m). c\] and thus define
 $n$ coordinate  forms $(\phi_c)_{c\in C}$ in $\hat M$. We check that
 \[\phi_c(m)=\phi_e(c^{-1}.m)\] so the $\bK$-linear map
  \begin{equation*}
\xymatrix@R-2pc{
  \check M     \ar@{->}[r]  & \hat M \\
  \phi     \ar@{|->}[r]  & \phi_e \\
}
\end{equation*}
  is an isomorphism of $\bK$-vector spaces. For
  every   $\psi$  in $\hat M$ we write
  $\psi^C$ for the corresponding element in $\check M$. So
  \[\psi^C(m)=\sum_{c\, \in \,C}\psi(c^{-1}.m)c.\]

  We now let  $L$ and $M$ 
be two  finitely generated $\bK[C]$-module.
We assume that $M$ is free.
So there exists
 a $\bK$-vector space $W$ such that \[M=\bigoplus\limits_{c\, \in \, C}c.W\]
 as a $\bK$-vector space.
 Let $t : L\times L\rightarrow M$ be a $C$-equivariant
 $\bK$-bilinear map. There are $n$  maps
 $(t_c)_{c\, \in \, C}$ such that
$t_c : L\times L\rightarrow W$
 is $\bK$-bilinear for every $c$ in $C$ and 
 for  $x$
 and $y$ in $L$ we have 
 \[t(x,y)=\sum_{c\, \in \, C}c.t_c(x,y).\]
We check that \[t_c(x,y)=t_e(c^{-1}.x, c^{-1}.y)\]
for every  $c\in C$  and $x$, $y$ in
$L$.
The map
  \begin{equation*}
\xymatrix@R-2pc{
   \Bil_C(L,M)    \ar@{->}[r]  & \Bil_\bK(L,W)\\
t  \ar@{|->}[r] & t_e
}
  \end{equation*}
  \smallskip
is thus an isomorphism
    between the $\bK$-vector space
    $\Bil_C(L,M)$ of $C$-equivariant
 \hbox{$\bK$-bilinear} maps from $L\times L$ to $M$, and
 the 
 space $\Bil_\bK(L,W)$ of 
 $\bK$-bilinear maps from $L\times L$ to $W$.
 For every $u$
 in  $\Bil_\bK(L,W)$ we write $u^C$ the corresponding
 map in $\Bil_C(L,M)$. So
 \[u^C(x,y)=\sum_{c\, \in \, C}c.u(c^{-1}.x,c^{-1}.y).\]

Let $\alpha_1$ and $\alpha_2$ in $\check L$.
Let
$(\alpha_{1,c})_{c\in C}$
be the  $n$ forms
in  $\hat L$ such that 
\[\alpha_1 (l) =\sum_{c\,\in\, C} \alpha_{1,c}(l).c\]
for every $l$ in $L$.
We similarly define $n$ forms
$(\alpha_{2,c})_{c\in C}$
in  $\hat L$.
Let $w\in W$ and let  $t=\gamma_{w,\alpha_1,\alpha_2}$.
 Then for $l_1$ and $l_2$ in $L$
 we have
 \[t(l_1,l_2)=\left( \alpha_1(l_1)\diamond \alpha_2(l_2)\right).w=
 \sum_{c\in C}\alpha_{1,c}(l_1)\alpha_{2,c}(l_2)c.w.
 \]
 We deduce
 \[t_e(l_1,l_2)=\alpha_{1,e}(l_1)\alpha_{2,e}(l_2)w
\text{\,\,  so \,\, } t_e=\pi_{w,\alpha_{1,e},\alpha_{2,e}}.\]
Equivalently, if $\beta_1$
and $\beta_2$ are in $\hat L$ and $w$ is in $W$ we have
\[\pi_{w,\beta_1,\beta_2}^C =\gamma_{w,\beta_1^C, \beta_2^C}.\]

We deduce that
if $L$ and  $M$ are  $\bK[C]$-modules
with  $M$  free, and if
$t : L\times L\rightarrow M$
is a $C$-equivariant bilinear map, then
every decomposition of $t_e$ as a sum
of $k$ pure  $\bK$-bilinear maps
results in a decomposition of
$t$ as a sum of $k$
pure $C$-equivariant  $\bK$-bilinear maps. So
\begin{equation*}R_{\bK,C}(t)\leqslant R_\bK(t_e)\leqslant \rank (M)\times (\dim_\bK (L))^2.\end{equation*}
And in case $t$ is symmetric
\begin{equation}\label{eq:boundc}
  S_{\bK,C} (t)\leqslant S_\bK(t_e)\leqslant  \rank (M)\times
\dim_\bK(L)\times 
(3\dim_\bK(L)-1)/2.\end{equation}

\section{The  complexity  of multiplication in finite fields}\label{sec:cff}

Let $\bK$ be a finite field with $q$ elements
and let $\bL$ be
a degree $n\geqslant 1$  field extension of $\bK$. 
The multiplication
map $\times_\bL : \bL\times \bL\rightarrow \bL$ is $\bK$-bilinear
and symmetric.
Its  bilinear complexity
$R_\bK(\times_\bL)$ is usually denoted $\mu_{q}(n)$
and its symmetric  bilinear complexity
$S_\bK(\times_\bL)$ is denoted $\mu_{q}^{ sym}(n)$.
It is known that $\mu_{q}^{ sym}(n)\geqslant 2n-1$. See \cite[Lemma 1.9.]{RAN}
for example. Lagrange 
interpolation
shows that  $\mu_{q}^{ sym}(n)=  2n-1$
when  $q\geqslant 2n-2$.
Chudnovsky and Chudnovsky have proved \cite{CC}  linear upper bounds 
for these bilinear complexities using interpolation
on algebraic curves.
Their method  
has been extensively studied and improved,
notably
by Shparlinski, Tsfasmann, Vladut \cite{STV},
Shokrollahi \cite{SH}, Ballet  and Rolland \cite{BR,BAL},
Chaumine \cite{CHA},  Randriambololona \cite{RAN}
and others, 
achieving sharper and sharper upper bounds for the bilinear 
complexity of multiplication in finite  extensions of finite fields.
See \cite{related} for a recent survey. We will use the following theorem.
\begin{theorem}[Chudnovsky (1987), Shparlinski, Tsfasmann and Vladut (1992), Ballet (1999)]\label{th:chud}
  There exists an effective absolute constant $\cQ$
  such that $\mu_q^{sym}(n)\leqslant \cQ n$ for all
  $n\geqslant 1$ and all prime power $q$.
  \end{theorem}

\section{The algebra $\bK[x]/x^n$}\label{sec:trunc}

 Let $\bK$ be a field with $q$
 elements. Let $n\geqslant 1$  be an integer.
 Let $\bL$ be a degree $2n-1$ field
 extension of $\bK$.
 Let $\bK[x]_{n-1}$
 be the $\bK$-vector space of polynomials with degree
 $\leqslant n-1$. The multiplication map $\bK[x]_{n-1}\times
 \bK[x]_{n-1}\rightarrow \bK[x]_{2n-2}$ is a restriction
 of the multiplication map $\bL\times \bL\rightarrow \bL$.
 And the multiplication map $\bK[x]/x^n\times
 \bK[x]/x^n\rightarrow \bK[x]/x^n$ is a restriction
 of $\bK[x]_{n-1}\times
 \bK[x]_{n-1}\rightarrow \bK[x]_{2n-2}$. So the symmetric
bilinear  complexity of multiplication in
 the quotient $\bK[x]/x^n$ is bounded from
 above by $\mu_q^{sym}(2n-1)$. So
 \begin{equation}\label{eq:rand}
   S_\bK(\times : \bK[x]/x^n\times
 \bK[x]/x^n\rightarrow \bK[x]/x^n)\leqslant \cQ n\end{equation}
 for some effective absolute constant $\cQ$.

 In case  $q\geqslant 2n-2$,
Lagrange  interpolation shows that the symmetric bilinear complexity of
$\bK[x]_{n-1}\times
 \bK[x]_{n-1}\rightarrow \bK[x]_{2n-2}$
 is $\leqslant 2n-1$. So the symmetric
bilinear  complexity of multiplication
in $\bK[x]/x^n$ is  $\leqslant 2n-1$ in that case.
In the other direction, Winograd has proved in \cite{WIN}
that this complexity is always $\geqslant 2n-1$. 
 More precise, more general  and stronger 
 statements can be found in \cite{BCS, BR, RAN} and
 \cite[Section 2]{related}.

  \section{The equivariant complexity  of multiplication in finite fields}\label{sec:ecff}

 Let $\bK$ be a finite field with $q$ elements
and let $\bL$ be
a degree $n\geqslant 1$  field extension of $\bK$. Let $C$
be the Galois group of $\bL/\bK$.
Then  $\bL$ is   a free $\bK[C]$-module of rank one.
We denote \[\times_{\bK, n}  : \bL\times \bL\rightarrow \bL\]
the multiplication map in $\bL$.
This is  a $C$-equivariant $\bK$-bilinear map.
We  define $\nu_q(n)$ to be the $C$-equivariant complexity
of $\times_{\bK,n}$ over $\bK$. We similarly define $\nu_q^{sym}(n)$ to 
 be the $C$-equivariant symmetric complexity
  of $\times_{\bK,n}$ over $\bK$.

 \bigskip

 The equivariant complexity
 $\nu_q^{sym}(n)$ controls the computational difficulty of
 multiplying two elements in $\bL$ given
 by their coordinates in a normal basis.
Indeed assume that $\nu_q^{sym}(n)=\sigma$.
 There exist two $\bK[C]$-linear maps
\[\top  : \bL\rightarrow \left(\bK[C]\right)^{\sigma} \text{ \,\, and \,\,}
 \bot  : \left(\bK[C]\right)^{\sigma} \rightarrow  \bL\]
 such that \begin{equation}\label{eq:tensor1}
   l_1\times l_2 =\bot (\top (l_1)\diamond_\sigma \top (l_2))\end{equation}
 for any $l_1$, $l_2$ in $\bL$.
 We note that
 $\top$ is
 a
 linear map
 between two
 free
 $\bK[C]$-modules
 of respective ranks
$1$ and $\sigma$.
Once chosen a basis
of $\bL$ 
we can describe $\top$
by  a $\sigma \times 1$ matrix with coefficients
in $\bK[C]$. Giving a basis
of $\bL$
as a $\bK[C]$-module boils down to
choosing  a normal
basis of $\bL/\bK$.
Similarly $\bot$
can be described by a $1\times \sigma$
matrix with coefficients
in $\bK[C]$.
So using Equation (\ref{eq:tensor1}) we  
compute the  product of two elements
in $\bL$ given by their coordinates in a
given normal basis in three steps:

 \begin{enumerate}[label=\arabic*.]
 \item Apply $\top $ to each element.
 \item Multiply the two elements thus obtained
   in $\bK[C]^\sigma$ using the $\diamond_\sigma$ law.
 \item Apply $\bot $
   to the result.
   \end{enumerate}

 The first step requires twice $\sigma $ multiplications
 in $\bK [C]$. We note that multiplication
 in $\bK [C]$ is the standard  convolution product. 
 The second step is a $\diamond_\sigma$  product between
 two  vectors in $\bK[C]^\sigma$. 
 The third step requires  $\sigma$ multiplications
 in $\bK [C]$.
  The  only bilinear step is the second one. All the multiplications
 in the first and third steps 
 involve  a variable and a constant.
 The total cost (omitting  additions) is $3\sigma $ convolution
 products between vectors of  length $n$ and
 $\sigma n$ multiplications in $\bK$. According to
 work by Sch\"{o}nhage and Strassen \cite{strassen}
 and Cantor and Kaltofen \cite{kal},
 convolution
 products of length $n\geqslant 2$ over an arbitrary commutative ring
 can be computed at the expense of
$O(n\log (n) |\log (\log (n))|)$ operations 
 in this ring. See also \cite[Theorem 8.23]{gathen}.
 Note in particular that it  is not necessary to have  $n$-th roots
 of unity in the base ring  in order to compute convolution
 products efficiently.

  \section{Semi-simple algebras over  finite fields}\label{sec:semi}

 Let $\bK$ be a finite field with $q$ elements. Let $n_1 \geqslant 1$ be an integer.
Let $\bL$ be
a degree $n_1$  extension of $\bK$. Let $F_q : \bL \rightarrow \bL$ be the Frobenius
automorphism of $\bL/\bK$.
Let $n_2 \geqslant 1$ be an integer. Set $\bM = \bL^{n_2}$.
This is  a semisimple $\bK$-algebra
of degree $n=n_1n_2$. 
We define an automorphism of $\bM$ over $\bK$ by sending
$(x_0, x_1, \ldots, x_{n_2-1})$
onto $(x_1, x_2, \ldots, x_{n_2-1}, F_q(x_0))$. We call $C$
the group generated by this automorphism. This is a cyclic
group of order $n$. And $\bM$ is a free $\bK[C]$-module of
rank $1$.
We let
\[\times_{\bK, n_1, n_2} : \bM \times \bM \rightarrow \bM\]
be the multiplication map in $\bM$.
This is a symmetric $C$-equivariant $\bK$-bilinear map.
We denote $\nu_q^{sym}(n_1,n_2)$ its symmetric
equivariant complexity.

 \section{Extension of scalars I}\label{sec:ext1}

 Let $\bK$ be a commutative field. Let $V$ and $W$ be two
 finite dimensional $\bK$-vector spaces. Let
 $t : V\times V\rightarrow W$ be a symmetric $\bK$-bilinear map.
Let  $S_{\bK}(t)$ be  the symmetric   bilinear complexity
of $t$.
 Let $\bL$ be a finite field  extension of $\bK$. We
 set $V_\bL = V\otimes_\bK \bL$,
 $W_\bL = W\otimes_\bK \bL$, $t_\bL = t\otimes_\bK\bL$.
Let  $S_\bL(t_\bL)$ be  the  symmetric bilinear complexity of
  $t_\bL$ as an 
 $\bL$-bilinear map. We have \begin{equation}\label{eq:extscal}
   S_\bL(t_\bL)\leqslant S_\bK(t).\end{equation}

We denote  by $S_\bK(\times_\bL)$ the symmetric  $\bK$-bilinear complexity
 of  $\times_\bL : \bL\times \bL\rightarrow
 \bL$,  the multiplication map in $\bL$ seen as a $\bK$-bilinear map.  According to  \cite[Lemma 1.10]{RAN} 
 \begin{equation}\label{eq:descscal}
   S_\bK(t)\leqslant S_\bL(t_\bL)\times S_\bK(\times_\bL).\end{equation}


 \section{Extension of scalars II}\label{sec:ext2}

 We now study the effect
 of extension of scalars on
 equivariant bilinear maps.
 The main motivation for extending  scalars is to increase the number
 of rational points  in the context of the geometric methods presented
 in Section \ref{sec:geo}.
Let $C$ be a finite  group of order $n$.
Let $\bK$ be a commutative field.
Let $L$ and $M$ be two finitely generated
$\bK[C]$-modules. Let $t : L\times L \rightarrow M$
be a symmetric
$C$-equivariant $\bK$-bilinear  form.
We denote $S_{\bK, C}(t)$ the symmetric equivariant  complexity
of $t$.

 Let $\bL$ be a finite field  extension of $\bK$. We
 set $L_\bL = L\otimes_\bK \bL$,
 $M_\bL = M\otimes_\bK \bL$, $t_\bL = t\otimes_\bK\bL$.
 We call $S_{\bL, C}(t_\bL)$ the  symmetric equivariant complexity of
  $t_\bL$. We have \begin{equation}\label{eq:extscalC}
   S_{\bL, C}(t_\bL)\leqslant S_{\bK, C}(t).\end{equation}

Let  $S_\bK(\times_\bL)$ be  the symmetric  bilinear  complexity
 of  $\times_\bL : \bL\times \bL\rightarrow
 \bL$, as a $\bK$-bilinear map. Then
 \begin{equation}\label{eq:descscalC}
   S_{\bK,C}(t)\leqslant S_{\bL,C}(t_\bL)\times S_\bK(\times_\bL).\end{equation}


Assume $\bK$ is a finite field with $q$ elements.
Let  $\bL$ be a degree $n$ field extension of
$\bK$. Let $C$ be the Galois group of $\bL/\bK$.
The multiplication
$\times_{\bK,n} : \bL\times \bL \rightarrow \bL$
is $\bK$-bilinear symmetric  and $C$-equivariant.
Let  $\bK'$ be a degree $m$ field extension of
$\bK$. We  tensor product the multiplication
$\times_{\bK,n}$ 
by $\bK'$ over $\bK$. The resulting
$C$-equivariant $\bK'$-bilinear map is isomorphic
to $\times_{\bK', n_1, n_2}$  as defined in Section \ref{sec:semi},
with
\[n_2=\gcd (n,m) \text{\,\,\, and \,\,\,} n_1=n/n_2.\]
Equations   \ref{eq:extscalC}
and
\ref{eq:descscalC} thus imply
\begin{equation*}
\nu_{q^m}^{sym}(n_1,n_2)\leqslant \nu_{q}^{sym} (n)
\end{equation*}
and
\begin{equation}\label{eq:split2}
\nu_q^{sym}(n)\leqslant \nu_{q^m}^{sym} (n_1,n_2)\times \mu_q^{sym }(m).
\end{equation}

 \section{A geometric construction}\label{sec:geo}

Let $\bK$ be a finite field with $q$ elements. We call
$p$ the characteristic of $\bK$. Let $Y$
be a smooth
absolutely  integral projective 
curve over $\bK$. Let $C$ be a cyclic group
of $\bK$-automorphisms of $Y$. 
We call $n$ the cardinality of $C$.
We assume that $n\geqslant 2$  and $p$ does not divide $n$.
We call  $X$ the quotient $Y/C$.
This is
a smooth
absolutely  integral projective 
curve over $\bK$.
We call
$\rho : Y \rightarrow X$ the quotient map.
Let $\rgot$ be an effective divisor on $Y$. We assume
that $\rgot$ and $c.\rgot$ are disjoint for every $c$ in $C$.
We set \[R=\sum_{c\in C} c.\rgot\] and call $\bK[R]$ the residue
ring at $R$. We identify  the ring $K[\rgot]$ with the subring
of $\bK[R]$ consisting of functions
vanishing at $c.\rgot$ for every $c$ in $C$ different from
$e$.
As a $\bK$-vector space
\[\bK[R]=\bigoplus_{c\in C}\, c.\bK[\rgot ].\]
 So $\bK[R]$ is a free
$\bK[C]$-module.
The
multiplication map
$\bK[\rgot ]\times \bK[\rgot ]\rightarrow \bK[\rgot ]$ is
$\bK$-bilinear and
symmetric.
We denote $\sigma$ its symmetric
bilinear complexity. According to Equation~(\ref{eq:boundc}) this is an upper
bound for the $C$-equivariant symmetric complexity
of $\bK[R]\times \bK[R]\rightarrow \bK[R]$. So there exist
two $\bK[C]$-linear maps 
 \[\top  : \bK[R]\rightarrow \bK[C]^{\sigma} \text{ \,\, and \,\,}
\bot  : \bK[C]^{\sigma} \rightarrow  \bK[R]\]
 such that \[l_1 \times l_2 =\bot (\top (l_1)\diamond_\sigma  \top (l_2))\]
 for $l_1, l_2  \in \bK[R]$.

We denote  $X(\bK)$ the set of $\bK$-points  on $X$.
 Let $a\in X(\bK)$ such  that $\rho$ is not ramified
 above $a$. Let $n_1$ be the inertial degree of $\rho$ at $a$.
 Let $n_2=n/n_1$. The fiber $B=\rho^{-1}(a)$ is
 a reduced  $\bK$-scheme consisting of $n_2$
 irreducible components, each of degree $n_1$ above $a$.
We call $\bM$ the residue ring $\bK[B]$ of $B$. 
This is a free $\bK [C]$-module of rank $1$.
As  a $\bK$-bilinear symmetric
$C$-equivariant map, the multiplication map in $\bM$ is
isomorphic to the map $\times_{\bK, n_1, n_2}$
introduced in Section \ref{sec:semi}. Its symmetric equivariant complexity
is thus $\nu_q^{sym}(n_1,n_2)$.

Let $D$ be a divisor on $X/\bK$. 
We call $E=\rho^{-1}(D)$
the pullback of $D$ on $Y$. 
Let $\epsilon$ be a local equation of $D$ in a neighborhood
of $a$ and $\rho(R)$. Seen as a function on $Y$ this is a local
 equation of $E$ in a neighborhood
of $B$ and $R$.  
Let
  \begin{equation*}
\xymatrix@R-2pc{
  e_B   & \relax : &     H^0(Y,\cO_Y(E))   \ar@{->}[r]  & \bM \\
  &&f   \ar@{|->}[r] & ({f}\times {\epsilon}) \bmod B}
  \end{equation*}
  be the evaluation map at $B$.  We similarly define
\begin{equation*}
\xymatrix@R-2pc{
  e_B^2   & \relax : &     H^0(Y,\cO_Y(2E))   \ar@{->}[r]  & \bM\\
  &&f   \ar@{|->}[r] & ({f}\times {\epsilon ^2}) \bmod B}
  \end{equation*}
  These maps are morphisms
  of $\bK [C]$-modules.
  For $f_1$ and $f_2$ in $H^0(Y,\cO_Y(E))$ we have
  \[e_B^2(f_1\times f_2)=e_B(f_1)\times e_B(f_2).\]
We assume that
$e_B$ is surjective. 
Since  $p$ does not divide
$n$, the ring $\bK[C]$ is semi-simple. So the kernel
of $e_B$ is a direct factor. We deduce that
$e_B$ has a right  inverse
\[e_B^\star : \bM\rightarrow H^0(Y,\cO_Y(E))\]
 which is  $\bK[C]$-linear. 
 Let
  \begin{equation*}
\xymatrix@R-2pc{
  e_R   & \relax : &     H^0(Y,\cO_Y(E))   \ar@{->}[r]  & \bK [R] \\
  &&f   \ar@{|->}[r] & ({f}\times {\epsilon }) \bmod R}
  \end{equation*}
and 
 \begin{equation*}
\xymatrix@R-2pc{
  e_R^2   & \relax : &     H^0(Y,\cO_Y(2E))   \ar@{->}[r]  & \bK [R] \\
  &&f   \ar@{|->}[r] & ({f}\times {\epsilon ^2}) \bmod R}
  \end{equation*}
 be the evaluation maps at $R$.
 These are  $\bK[C]$-linear maps.
We assume that $e_R^2$
is injective. 
Since  the ring $\bK[C]$ is semi-simple, the image
of $e_R^2$ is a direct factor of $\bK[R]$. 
We deduce the existence
of a  left inverse \[e_R^\star :  \bK [R] \rightarrow  H^0(Y,\cO_Y(2E))\]
to the evaluation map $e_R^2$.
Let $s_1$ and $s_2$ be two funtions in $H^0(Y,\cO_Y(E))$, representing the two
elements \[e_B(s_1)= (s_1\times \epsilon )
\bmod B \text{ \,\, and \,\, } e_B(s_2)=( s_2\times \epsilon ) \bmod B\]
in $\bM$.
The product $s_3=s_1s_2$ belongs to $H^0(Y,\cO_Y(2E))$
and \[e_R^2(s_3)=e_R(s_1)\times e_R(s_2) \in \bK[R].\] So
\[s_3=e_R^\star(e_R(s_1)\times e_R(s_2))
  =e_R^\star(\bot (\top (e_R(s_1))\diamond_\sigma \top (e_R(s_2))))\] and
the \hbox{$\bK$-bilinear} map
\begin{equation*}
  e_B^2 \circ  e_R^\star \circ \bot \circ \diamond_\sigma \circ (\top
  \times \top )  \circ (e_R\times e_R) \circ (e_B^{\star}\times
e_B^{\star}): \bM\times \bM \rightarrow \bM\end{equation*}
is the multiplication map in $\bM$.
We observe that \[\top \circ e_R\circ e_B^{\star}  : \bM \rightarrow \bK[C]^\sigma\]
and
\[e_B^2\circ  e_R^\star \circ \bot : \bK [C]^\sigma\rightarrow \bM\]
are  $\bK[C]$-linear maps.
We deduce that
\begin{equation*}
  \nu_q^{sym}(n_1,n_2)\leqslant \sigma.\end{equation*}

 \section{A general bound}\label{sec:heu}

 Let $\bK$ be a finite field with $q$ elements.
 Let $\bKb$ be an algebraic closure of $\bK$.
 Let $n\geqslant 2$ be a prime to $q$ integer.
 Let $n_1$ and $n_2$ be two positive integers such that $n=n_1n_2$.
 We would like  to instantiate the construction in Section~\ref{sec:geo}
 so as to obtain a sharp bound for the
 equivariant symmetric complexity $\nu_q^{sym}(n_1,n_2)$ of multiplication
 in the  degree $n$ algebra over $\bK$ defined in Section \ref{sec:semi}.
 Field extensions correspond to the case $n_2=1$.
 We let  $X$ be  
 a smooth absolutely integral curve over $\bK$ such that 
 \begin{equation}\label{eq:exa}
    X(\bK)\not = \emptyset.
 \end{equation}
 Let $\mu_n$ be a primitive $n$-th root of unity in $\bKb$.
 Let $\bK (\mu_n)$ be the field generated by
 $\mu_n$ over $\bK$. We assume that the Jacobian $J_X$ has a point
 \begin{equation}\label{eq:ext}s \in J_X(\bK (\mu_n)) \text{\, of order
     \,} n, \text{\,  such that \,} F_q(s)=qs\end{equation}
 where  $F_q$ is the Frobenius  of  $J_X/\bK$. A sufficient condition
 for such  an $s$ to exist is that  the characteristic polynomial $\chi (t)$  of
  $F_q$  has a root in $\bZ_\ell$
 congruent to $q$ modulo $n$, for every  prime $\ell$
 dividing $n$. 
 This is granted if $n$
 divides the cardinality $\chi(1)$
 of $J_X(\bK)$, and $1$ is a simple root
 of $\chi$  modulo $\ell$ for every prime $\ell$
 dividing $n$. That is \begin{equation}\label{eq:cnum}
   \chi(1)=0\bmod n \text{ \,\, and \,\,} \gcd (\chi'(1),n)=1\end{equation}
   where $\chi'$ is the derivative of the polynomial $\chi$.
We look for a curve $X$ with smallest possible
genus satisfying these conditions.
Condition
(\ref{eq:cnum})
cannot hold if $n>(1+\sqrt q)^{2g}$. On the
other hand we  heuristically expect  to find a
curve
$X$ with genus $g_X$ equal to $g$  and satisfying condition
(\ref{eq:cnum}) provided
\begin{equation*}
  g\gg \log_q n.\end{equation*}
Conditions  (\ref{eq:exa})
and (\ref{eq:ext}) and Kummer theory
imply  the existence
of a curve $Y$ over $\bK$
and an unramified  Galois cover $\rho : Y\rightarrow X$
with cyclic Galois group of order $n$. We can even force
a $\bK$-point on $X$ to split completely in $Y$.
We take $a\in X(\bK)$, $B=\rho^{-1}(a)$, $n_1$, $n_2$,
$\rgot$, 
$R=\sum_{c\in C}c.\rgot$,  $D$, $E=\rho^{-1}(D)$,  $e_B$
 and $e_R$ as in Section~\ref{sec:geo}. 
The  condition 
 \begin{equation}\label{eq:injsur}
e_B \text{\,\, is surjective and \,\, }
e_R^2 \text{\,\, is injective}\end{equation}
is granted if 
\[\deg (E-B)> 2g_Y-2  \text{\,\, and \,\,}  \deg (2E-R)<0\]
or equivalently 
\begin{equation*}
  \deg D \geqslant 2g_X \text{\, and \,}
  \deg \rgot \geqslant 2\deg D+1.
\end{equation*}

This last condition is easy to check
but a bit  restrictive.
A more delicate sufficient condition  for (\ref{eq:injsur}) is
\begin{equation*}
  E-B \text{\,\, is non-special and \,} \dim H^0(Y,\cO_Y(2E-R))=0.
\end{equation*}

Remind that a divisor $D$ on a curve $X$ of genus $g_X$
      is said to be {\bf non-special} if the dimension  of $H^0(X,\cO_X(D))$
    is $\deg D -g_X+1$. Otherwise $D$ is said to be {\bf special}.
We summarize the above discussion in the theorem bellow.

\begin{theorem}\label{th:analy}
  Let $\bK$ be a finite field with $q$ elements. Let
  $n\geqslant 2$ be a prime to $q$
  integer. Let $\rho : Y\rightarrow X$
  be an unramified  Galois  cover between two smooth
  absolutely
  integral curves
  over $\bK$. We assume that the Galois group $C$
  of $\rho$ is cyclic of order $n$. 
  Let $a\in X(\bK)$.  
   Let $n_1$ be the inertial degree of $\rho$ at $a$. Let   $B=\rho^{-1}(a)$ be the fiber
  of $\rho$ above $a$. 
          Let $n_2=n/n_1$.
  Let $\rgot$ be an effective divisor on $Y$ such that
  $\rgot$ and $c.\rgot$ are disjoint for every $c$ in $C$.
Let $R=\sum_{c\in C}c.\rgot$.
  Let $D$ be a divisor on $X/\bK$.  Let  $E=\rho^{-1}(D)$.
  We assume that
 \begin{equation}\label{eq:numcrit}
  \deg D \geqslant 2g_X \text{\, and \,}
  \deg \rgot \geqslant 2\deg D+1.
\end{equation}
or
\begin{equation}\label{eq:nonspecial}
  E-B \text{\,\, is non-special and \,\, } \dim H^0(Y,\cO_Y(2E-R))=0.
\end{equation}
  where $g_X$ is the genus of $X$. Then $\nu_q^{sym}(n_1,n_2)\leqslant \sigma$
  where \[\sigma = S_\bK(\times : \bK[\rgot ]\times \bK[\rgot ]
  \rightarrow \bK[\rgot ])\] is the symmetric bilinear complexity
  of multiplication in the residue ring of $\rgot $.\end{theorem}

If  $\rgot $ is $\deg \rgot $ times a point
in $Y(\bK)$, the symmetric
bilinear   complexity $\sigma$ 
of $\bK[\rgot ]\simeq \bK[x]/x^{\deg \rgot}$ is linear in the degree of $\rgot$
according to Equation~(\ref{eq:rand}).
If $\rgot $ is reduced and irreducible then $\sigma$   is linear in the degree of $\rgot$
according to Theorem \ref{th:chud}. If $\rgot $ is  a sum of
$\deg \rgot $ pairwise distinct $\bK$-rational points then
$\sigma =\deg \rgot $.

\section{Non-special divisors}\label{sec:NS}

In order
to verify
Condition (\ref{eq:nonspecial})
in Theorem \ref{th:analy},
we need a simple criterion for a divisor to be non-special
in this context.
Let $\bK$ be a 
field with characteristic $p$. Let $n\geqslant 2$
be a prime to $p$ integer.
Let $\bKb$ be an algebraic closure of $\bK$.
 Let $\mu_n$ be a primitive $n$-th root of unity in $\bKb$.
Let $X$ and $Y$ be two smooth absolutely
integral curves over $\bK$.
We call $g_X$ the genus of $X$ and
$g_Y$ the genus of $Y$. Let $\rho : Y\rightarrow X$ be a Galois unramified
cover with cyclic
Galois group $C$ of order $n$.
Let $\hat\rho : J_X\rightarrow J_Y$
be the induced map on Jacobian varieties.
The kernel of $\hat\rho$ is a  finite  group scheme 
of degree  $n$. There is a  pairing
\[e_\rho : \Ker \hat\rho \times C \rightarrow \mmu_n\]
where $\mmu_n$ is the group scheme  of $n$-th roots of unity.
If $\gamma$
is a divisor class in the kernel of $\hat \rho$ and $c\in C$,
we let $\Gamma$ be a divisor in $\gamma$ and $G$ a function
on $Y$ with divisor $\rho^{-1}(\Gamma)$. We set \[e_\rho(\gamma , c)=
  \frac{G\circ c}{G}.\]
This is a non-degenerate pairing. As a  consequence  $\Ker \hat\rho$
is isomorphic to $\mmu_n$.

Let  $D$ be a  divisor on $X$ with degree $g_X-1$.
Let  $E$ be the pullback of $D$ on $Y$.
The degree of $E$ is $n(g_X-1)=g_Y-1$.
If  $E$ is special then the $\bK[C]$-module
$H^0(Y,\cO_Y(E))$
is non-zero. Since $p$ is prime to $n$,
there exists an eigenvector $\varphi$
for the action of $C$ on 
$H^0(Y,\cO_Y(E))\otimes \bK (\mu_n)$. 
 Let $N$ be the effective divisor on $Y\otimes \bK (\mu_n) $
such that the principal divisor $(\varphi)$ is $N-E$. Then $N$ is the pullback of an effective
divisor $M$ on $X\otimes \bK (\mu_n)$. And $M-D$ is in the the kernel of
$\hat\rho$.

To summarize, if  $D$ is
a  divisor on $X$ with degree $g_X-1$, then
the pullback $E=\rho^{-1}(D)$
has degree $g_Y-1$. And $E$  is special (its class  is effective) if and
only if there exists a degree $g_X-1$ effective
divisor $M$ on $X\otimes \bK (\mu_n)$ such that $D-M$ is  in the kernel of
$\hat\rho$.

\section{Elliptic curves}\label{sec:ec}

In this section we adapt the general method of Section~\ref{sec:geo}
to the special case of elliptic curves. The main reason
for this restriction is that we have a good control
on the group of rational points
on an elliptic curve. 
Restricting to elliptic curves is not optimal 
 but it  enables
us  to prove
such an asymptotic statement as Theorem \ref{th:asymp}.

Let $\bK$ be a finite field with cardinality $q$ and
characteristic $p$. Let $n\geqslant 2$ be  an integer.
Let $Y$ be 
an elliptic curve over $\bK$. We assume
that $Y$ has a $\bK$-point $t$
of order $n$. Let $C$ be the group generated by $t$.
Let  $X$ be  the quotient of $Y$ by  $C$.
Let  $\rho : Y\rightarrow X$ be the quotient
isogeny. Let  $a$ in $X(\bK)$.
Let $n_1$ be the inertial degree of $\rho$  at $a$. Let $n_2=n/n_1$. The  fiber $B=\rho^{-1}(a)$
of   $\rho$ above   $a$  has
$n_2$ irreducible components, each of degree $n_1$ above $a$. We call
  $\bM$ the residue ring $\bK[B]$ of $B$. 
Let $v$ be a point
in  $X(\bK)$. 
We assume that $v-a$ is not in the kernel
  of
the dual isogeny
  $\hat\rho$.  
   We call $D$
   the degree $1$ divisor on $X$ consisting of the
   single point $v$ with multiplicity $1$.
   We call $E$ the divisor $\rho^{-1}(D)$.
      We let $u$ be a non-zero point
      in  $\rho(Y(\bK))$.
      We  assume that $u-2v$ is not 
in the kernel
  of $\hat\rho$.  
   We let $\rgot$ be the formal  sum
   of one point in  the fiber 
   of $\rho$ above $u$ plus  one point in the kernel of $\rho$.
   Let    $R$ be the  closure of $\rgot $
   under   the action of $C$. So $R$ is the 
   sum of the two split fibers of $\rho$ above the origin
   $o_X$ and $u$.
   We
let $\epsilon$ be
   a local equation of $D$ in a neighborhood
   of $o_X$, $a$ and $u$.
   In this setting the evaluation map $e_B : H^0(Y,\cO_Y(E))\rightarrow  \bM$
  is an isomorphism between two free  $\bK[C]$-modules of rank $1$.
  The evaluation map $e^2_R : H^0(Y,\cO_Y(2E))\rightarrow  \bK[R]$
  is an isomorphism between two free  $\bK[C]$-modules of rank $2$.
The
  symmetric bilinear complexity of  multiplication in
  $\bK[\rgot ]\sim \bK \times \bK$ is $2$.
  We deduce that $\nu_q^{sym}(n_1,n_2)$ is bounded by $2$.
  In the special case when $n_1=n\geqslant 2$ the residue ring $\bM$
  is a field and   $\nu_q^{sym}(n_1,n_2)=\nu_q^{sym}(n)$. The latter  cannot be equal to $1$ then
  because  $\bK^n$ is not a field. So  $\nu_q^{sym}(n)=2$ in that case.

       \begin{theorem}\label{th:elcons}
    Let $\bK$ be a field with cardinality $q$ and
    characteristic $p$. Let $n\geqslant 2$ be  an integer.
    Let $Y$ be an elliptic curve over $\bK$ having a $\bK$-point $t$
    of order $n$.
    Let  $X$ be  the quotient of $Y$ by the group $C$
generated by $t$. Let  $\rho : Y\rightarrow X$ be the quotient
isogeny. Let $a\in X(\bK)$.
Let $n_1$ be the inertial degree  of $\rho$ at $a$. Let $n_2=n/n_1$.
Let $u$ be a non-zero point
in $\rho(Y(\bK))$. Let $v\in X(\bK)$. We assume that
neither $v-a$ nor $u-2v$ are 
in the kernel   of the dual isogeny
$\hat\rho$.  Then $\nu_q^{sym}(n_1,n_2)\leqslant 2$.     In case
$n_2=1$  then $\nu_q^{sym}(n)=2$.
  \end{theorem}

\section{An asymptotic bound}\label{sec:gb}

Let $n\geqslant 2$ be an integer.
Using
Theorem \ref{th:elcons} we now  prove an asymptotic 
bound on  the  equivariant
complexity $\nu_q^{sym}(n)$ without any restriction on $q$ or $n$.
We let $\bK$ be a field with cardinality $q$ and
characteristic $p$. Let $n\geqslant 2$ be  an integer.
We first assume  that \begin{equation}\label{eq:condn1}
  n^2\leqslant 2\sqrt q\end{equation} and
  \begin{equation}\label{eq:condn2}
    q\geqslant 37.\end{equation} There are
two consecutive multiples of $n^2$ in the Hasse interval
$[q+1-2\sqrt q, q+1+2\sqrt q]$. At least one of them
is not congruent to $1$ modulo $p$. So there exists
an elliptic curve $Y$ over $\bK$ such that $Y(\bK)$ is
divisible by $n^2$. We deduce that $Y$ has a $\bK$-point
of order $n$. Indeed the group $Y(\bK)$ is isomorphic to
$(\bZ/m_1\bZ)\times (\bZ/m_1m_2\bZ)$ where $m_1$
and $m_2$ are positive integer. And $n^2$ divides
$\# Y(\bK)=m_1^2m_2$.  So $n^2$ divides $(m_1m_2)^2$.
So $n$ divides $m_1m_2$. So $Y(\bK)$, beeing isomorphic
to $(\bZ/m_1\bZ)\times (\bZ/m_1m_2\bZ)$,  has an element
$t$ of order $n$. Let $C$ be the group
generated by $t$.
We call $X$ the quotient of $Y$ by  $C$.
Let  $\rho : Y\rightarrow X$ be the quotient
isogeny. Let $P$ be a point in $X(\bK)$.
Let $\bKb$ be an algebraic closure of $\bK$.
Let
$Q$ be any point in $Y(\bKb)$ such that
$\rho(Q)=P$. We set $\kappa (P)=F_q(Q)-Q$
where $F_q$ is the Frobenius endomorphism of $Y/\bK$.
We thus define a morphism
  \begin{equation*}
\xymatrix@R-2pc{
  \kappa   & \relax : &    X(\bK)/\rho(Y(\bK))    \ar@{->}[r]  & \Ker \rho = \, C\\
  &&P   \ar@{|->}[r] & F_q(Q)-Q}
  \end{equation*}
  which is easily seen to be a bijection. 
Let $n_1$ and $n_2$ be two positive integers such that $n=n_1n_2$.
  There exists
 at least one 
  point $a$ in $X(\bK)$ such that
  $\kappa (a)=n_2t$.
  The  inertial degree of $\rho$ above $a$
  is $n_1$. 
  We call  $B=\rho^{-1}(a)$ the fiber of
  $\rho$ above   $a$.
It has 
$n_2$ irreducible components, each of degree $n_1$ above $a$. We call
  $\bM$ the residue ring  $\bK[B]$ of $B$.

  We need  a point  $v$ 
  in  $X(\bK)$
  such that
   $v-a$ is not in the kernel
  of $\hat\rho$.  There are at least $|X(\bK)|-n$ such points. So
  the existence
  of $v$ is granted provided
  \[|X(\bK)|-n\geqslant 1.\]
  The latter inequality follows from Conditions (\ref{eq:condn1}) and
   (\ref{eq:condn2}).
  We also need  a non-zero point  $u$ 
  in  $\rho(Y(\bK))$ 
  such that  $u-2v$ is not 
in the kernel
  of $\hat\rho$.  
  There are at least
 \[\frac{|X(\bK)|}{n}-n-1\]
 such points. So
  the existence
  of $u$ is granted provided
   \[\frac{|X(\bK)|}{n}-n\geqslant 2.\]
  The later inequality follows again  from Conditions (\ref{eq:condn1}) and
   (\ref{eq:condn2}). 
  Applying Theorem~\ref{th:elcons} we deduce the following.
  
  \begin{theorem}\label{th:large}
      Let $q$ be a prime power and $n\geqslant 2$
  an integer. Let $n_1$ and $n_2$ be two positive integers such that $n=n_1n_2$.
  If $q\geqslant 37$ and $n\leqslant \sqrt{2\sqrt q}$
  then $\nu_q^{sym}(n_1,n_2)\leqslant 2$.
\end{theorem}    

We now can bound $\nu_q(n)$ without any restriction
on $n$ and $q$. 
We let $m$ be the smallest integer such
that $m \geqslant 4\log_q n$ and $m \geqslant 6$.
We set $q'=q^m$ and check that $(q',n)$ satisfy
 Conditions
(\ref{eq:condn1}) and
(\ref{eq:condn2}). 
Using Theorem \ref{th:large} in conjunction  with Equation~(\ref{eq:split2}) and Theorem~\ref{th:chud}
we deduce the following theorems.

\begin{theorem}\label{th:asymp}
  Let $q$ be a prime power and $n\geqslant 2$
  an integer. Let
$m$ be the smallest integer such
that $m \geqslant 4\log_q n$ and $m \geqslant 6$.
  Then $\nu_q^{sym}(n)\leqslant 2\times \mu_q^{sym}(m)$.
\end{theorem}

\begin{theorem}\label{th:redu}
  There exists an absolute constant
  $\cQ$ such that the following is true.
  Let $q$ be a prime power and $n\geqslant 2$
  an integer. Then $\nu_q^{sym}(n)\leqslant \cQ \times \lceil \log_q n \rceil $.
\end{theorem}

The next theorem now   follows from Theorem \ref{th:redu}
and the existence of an  algorithm to compute 
  products in $\bK[x]/(x^n-1)$  at the expense of $O(n\log (n) |\log (\log (n))|)$ operations 
  in $\bK$.
  See \cite{strassen,kal}.

\begin{theorem}\label{th:thel}
  Let $\bK$ be a finite field of cardinality $q$. Let $\bL/\bK$ be
  an extension of degree $n\geqslant 2$. Let $\cB$ be a normal basis of $\bL/\bK$.
  There exists a straight line program that computes the coordinates
  in $\cB$ of the product of two
  elements in $\bL$ given by their coordinates in $\cB$ at the expense
  of \[\leqslant \cQ \times n\times \lceil \log_q(n) \rceil \times  \log (n)\times |\log(\log (n))|\] operations in $\bK$
  where $\cQ$ is an absolute constant.
\end{theorem}

  Compared to \cite[Theorem 4]{cel}
  we save
a $\log n$ factor on both the running time and the size of the model.
Theorem \ref{th:thel}
is also  more general since it applies  to any normal
basis and does not rely
on any ad hoc redundant  representation  as  in \cite{cel}.

\section{Bounding $\nu_q^{sym}(n)$}\label{sec:exp}

We explain
how to use
Theorems  \ref{th:elcons}
and \ref{th:analy}  to bound $\nu_q^{sym}(n)$ for
given $q$ and $n$.
If we plan to use an elliptic curve, we look for the smallest integer
$m$ such that the Hasse interval
\[[\lceil q^m+1-2q^{m/2}\rceil ,\lfloor q^m+1+2q^{m/2}\rfloor ] \] contains a multiple
of $n$. We then look for an elliptic curve over a field with $q^m$
elements satisfying the hypotheses  of  Theorem \ref{th:elcons}.
We pick random curves and compute their cardinality  using Schoof's algorithm and
its variants \cite{schoof},  until we find a  curve
with order divisible by $n$. We then check for the existence
of a point of order $n$.

If we want to use the general method of Section \ref{sec:geo},
we look for
the smallest integer  $g$ such that $(\sqrt q+1)^{2g}$ is reasonably larger than $n$.
We then pick random
curves of genus $g$ over a field with $q$ elements, until we find one
whose Jacobian has order divisible by $n$. We then check the hypotheses
of Theorem \ref{th:analy}.
We illustrate this method
with a few examples in the following sections.
We will see how to  verify the hypotheses
of Theorem \ref{th:analy} at  the least computational cost.
The knowledge of the zeta function suffices in many cases.

\section{The case $q=7$ and $n=5$}\label{sec:5}

Since $10$ belongs to the Hasse interval \[[ \lceil 7+1-2\sqrt 7 \rceil ,
\lfloor   7+1+2\sqrt 7 \rfloor ]= [3,13]\]
there is  an elliptic curve $E$
such that $E(\bK)\simeq \bZ/10\bZ$.
We can take $Y$
to be  the smooth projective model of  \[y^2=x^3+x+4.\] 
The point $t=(6,4)\in Y$ has order $5$. The quotient of $Y$
by the group $C$ generated by $t$ is the elliptic
curve $X$ with affine equation \[y^2 = x^3 + 3x + 4.\] Since
the kernel of the quotient by $C$ isogeny
$\rho : Y\rightarrow X$ is split, the kernel of the dual
isogeny $\hat\rho$ is isomorphic to $\mmu_5$. The only rational
point in it is the origin $o_X$ because $n$ is prime to $q-1$.
The image $\rho(Y(\bK))$ has order $2$.
The point \[a=(0,2)\in X(\bK)\] has order $5$. So it does
not belong to $\rho(Y(\bK))$. The fiber $B=\rho^{-1}(a)$ contains
no $\bK$-point. So it is irreducible.
We set \[u=(6,0)\in X(\bK)\] the unique $\bK$-rational point of
order $2$ on $X$. So $u$ belongs to $\rho (Y(\bK))$. We take
\[v=(0,5) \in X(\bK).\]
Since  $2v$ has order  $5$,  it must be different from $u$.
Since the only $\bK$-point in the kernel of
$\hat\rho$ is $o_X$ we easily check  that $v-a$ and $u-2v$ are not
in this kernel. Applying Theorem \ref{th:elcons}
we deduce that \[\nu^{sym}_7(5)=2.\]

\smallskip

The following  computer session implements
this calculation in 
SageMath (Version 9.4)  \cite{sagemath}.

\smallskip

\begin{verbatim}
sage: q=7;K=GF(q);n=5
....: Y=EllipticCurve([K(1),K(4)]);Y.order()
10
sage: t=Y(6,4);n*t
(0 : 1 : 0)
sage: rho = Y.isogeny(t);X = rho.codomain()
Elliptic Curve defined by y^2 = x^3 + 3*x + 4 
over Finite Field of size 7
sage: a=X(0,2);5*a
(0 : 1 : 0)
sage: u=X(6,0);v=X(0,5);u-2*v
(2 : 5 : 1)
\end{verbatim}

\section{The case  $q=11$ and $n=239$}\label{sec:239}

We  try the general method first. We then see what can be achieved
using  elliptic curves and extension of scalars.

\subsection{Using a genus $2$ curve}\label{sec:2391}

Let $X$ be the smooth
projective model  of the hyperelliptic curve with equation
\[y^2=x^5+x^3+2x^2+3.\] This is a genus $2$ curve.
The characteristic equation of the Frobenius $F_q$ of $X$
is
\[\chi(t)=t^4 + 7t^3 + 33t^2 + 77t + 121.\]
So $X$ has  $q+1+7=19$ points over $\bK$.
Its Jacobian has $\chi(1)=239=n$ points. This is a prime integer.
The factorization of $\chi(t)$ modulo $n$
is \[\chi(t) =  (t -11)(t -1)(t^2 + 19t + 11).\] So there is a point $s$
in $J_X[n]$ such that $F_q(s)=qs$.
Let $w_0\in X(\bK)$ be the  unique place  at infinity.
The class of $2w_0$ is the unique divisor class of degree $2$ on $X$
having 
positive projective dimension.
There exists a curve
$Y$ over $\bK$ and a Galois cover $\rho : Y \rightarrow X$ with
cyclic Galois group of order $n$ such that  the fiber
of $\rho$ above $w_0$ splits completely over $\bK$.
The kernel of $\hat\rho : J_X\rightarrow J_Y$ is the subgroup
generated by the class $s$. We observe that the class
of  $q$ generates
a subgroup of index $2$
in the multiplicative group $(\bZ/n\bZ)^*$. So Galois action 
on the non-zero classes in the kernel of $\hat\rho$ has two orbits.

Let $w_1$ and $w_2$ be two points in $X(\bK)$ having distinct
$x$-coordinates. The linear pencil of the divisor $w_1+w_2$ has
projective dimension zero, that is  \[H^0(X,\cO_X(w_1+w_2))=\bK.\]
Let $v_1$, $v_2$, $v_3$, $v_4$, $v_5$
be five  points in $X(\bK)$.
We assume that  $v_1$, $v_2$, $v_3$, $v_4$, $v_5$, $w_1$, and $w_2$
are  pairwise distinct.
Since the cardinality of $J_X(\bK)$ is odd,
the multiplication by two map is a bijection of it.  We deduce the existence
of five   effective degree two divisors  $D_1$, $D_2$, $D_3$, $D_4$, $D_5$
such that
$2(D_i-2w_0)$ is linearly equivalent to $w_1+w_2-v_i-w_0$ for
$1\leqslant i\leqslant 5$.
The divisor $2D_i-3w_0$ is linearly equivalent
to $w_1+w_2-v_i$. It is a non-special divisor.

Let $\xi$ be any non-zero
divisor class in the kernel of $\hat\rho$.
For each $1\leqslant i\leqslant 5$,
the divisor class $2D_i-3w_0-\xi$ is the class
 $w_1+w_2-v_i-\xi$. At most two  among these five  classes are effective.
Otherwise the class $w_1+w_2-\xi$ would have positive projective dimension.
So it would be  the class of $2w_0$. Then $\xi = w_1+w_2-2w_0$
would be  $\bK$-rational.
A contradiction.

Since there are only two Galois orbits on the non-zero classes
in $\Ker \rho$ we deduce that there exists a $v$
 among $v_1$,
 $v_2$, $v_3$, $v_4$, $v_5$ such that  $w_1+w_2-v-\xi$ is ineffective
 for all $\xi$ in this kernel.
We call
$D$ the effective degree two divisor  such that
$2(D-2w_0)$ is linearly equivalent to $w_1+w_2-v-w_0$.

Because $n$ is a prime integer, every  fiber of $\rho$
above a rational point of $X$ is either irreducible
or completely split. Since de genus of $Y$ is
\[g_Y=1+n(g_X-1)=1+n=240,\] the number of $\bK$-rational
points on it is bounded from above by $q+1+2g_Y\sqrt q <1604$.
So the number of split fibers is $\leqslant 1603/239<7$. So there are
at least $13$  points $(a_i)_{1\leqslant i\leqslant 13}$
in $X(\bK)$ with an irreducible
fiber above them.

At most two among  the
$(a_i)_{1\leqslant i\leqslant 13}$
make the class of $D-a_i$ effective. Otherwise $D$ would have
positive projective dimension. So it would be linearly
equivalent to $2w_0$. Then $w_1+w_2-v-w_0$ would be principal.
But $w_1+w_2$ has projective dimension $0$ and $v$ is distinct
from $w_1$ and $w_2$. A contradiction.

Let $\xi$ be any non-zero class
in the kernel of $\hat\rho$.
At most two among  the
$(a_i)_{1\leqslant i\leqslant 13}$
make $D-a_i-\xi$ effective. Otherwise $D-\xi$ would have
positive projective dimension. So it would be linearly
equivalent to $2w_0$. Then $\xi$ would be the
class of $D-2w_0$ and it would  therefore be  $\bK$-rational.
A contradiction.

Since  Galois action has two orbits on the non-zero
classes in the kernel of $\hat\rho$, we deduce that at least
$13-2-2\times 2 = 7$ rational points $a_i$ on $X$ have irreducible
fiber $\rho^{-1}(a_i)$ and make
$D-a_i+\xi$ non-special for every $\xi$ in $\Ker \hat\rho$.
We let $a$ be any of them.

We let  $\rgot$  be the three times any point on $Y$ above
$w_0$.
The ring $\bK[\rgot ]$ is isomorphic to $\bK[x]/x^3$. Since $\bK$ has $q\geqslant 4$
elements  the symmetric bilinear complexity of
multiplication in $\bK[\rgot ]$ is $5$. Applying Theorem \ref{th:analy}
we deduce \[\nu^{sym}_{11}(239)\leqslant 5.\]

\smallskip

The following  computer session implements
this calculation  in 
SageMath (Version 9.4)  \cite{sagemath}.

\smallskip

\begin{verbatim}
sage: q=11;K=GF(q);Kx.<x>=FunctionField(K);Kxy.<y>=Kx[];
KX.<y> = Kx.extension(y^2-x^5-x^3-K(2)*x^2-K(3)) 
g = KX.genus();LP=KX.L_polynomial();t=LP.parent().gen() 
sage: chi=LP(1/t)*t^(2*g)                                                
t^4 + 7*t^3 + 33*t^2 + 77*t + 121
sage: n=numerator(chi(1))   
239
sage: Fnt.<t> = PolynomialRing(GF(n));factor(Fnt(chi)) 
(t + 228)*(t + 238)*(t^2 + 19*t + 11)
\end{verbatim}


\subsection{Using an elliptic curve and extension of scalars}

We now try to bound $\nu^{sym}_{11}(239)$
using the method in Section \ref{sec:ec}. We let
$m$ be the smallest integer such that the Hasse interval
\[[\lceil q^m+1-2q^{m/2}\rceil ,\lfloor q^m+1+2q^{m/2}
  \rfloor ]\] contains a multiple
of $n$. For $m=1$ we find  the interval
$[9,15]$. For $m=2$ we find  the interval
$[111,133]$. For $m=3$ we find  the interval
$[1296,1368]$.  None of these three intervals
contain a multiple of $239$. So we must take $m\geqslant 4$.
The best  we can hope with this method is to prove
that \[\nu^{sym}_{11}(239)\leqslant 2\mu_q^{sym}(4).\] Since
$q\geqslant 6$ we have $\mu_q^{sym}(4) =  7$. So
$\nu^{sym}_{11}(239)\leqslant 14$. This is not as good as
the bound already  obtained in Section \ref{sec:2391}.

\section{The case  $q=13$ and $n=4639$}\label{sec:4639}

Let $X$ be the smooth projective
plane quartic  with homogeneous equation
\[y^3z + x^3y + 2xyz^2 + yz^3 + 11x^3z + 9x^2z^2 + 10xz^3.\]
This is a genus $3$ curve.
The characteristic equation of the Frobenius $F_q$ of $X$
is
\[\chi(t)= t^6 + 9t^5 + 51t^4 + 197t^3 + 663t^2 + 1521t + 2197.\]
So $X$ has  $q+1+9=23$ points over $\bK$.
Its Jacobian has $\chi(1)=4639=n$ points.
This is a prime integer.
The factorization of $\chi(t)$ modulo $n$
is \[\chi(t) =  (t -13)(t -1)(t+2195)(t+3726)(t^2+3380t+13).\] So there is a point $s$
in $J_X[n]$ such that $F_q(s)=qs$.
Let $w_0\in X(\bK)$ be the  point $(0:0:1)$.
There exists a curve
$Y$ over $\bK$ and a Galois cover $\rho : Y \rightarrow X$ with
cyclic Galois group of order $n$ such that  the fiber
of $\rho$ above $w_0$ splits completely over $\bK$.
The kernel of $\hat\rho : J_X\rightarrow J_Y$ is the subgroup
generated by the class $s$. We observe that the class
of  $q$ generates
the multiplicative group $(\bZ/n\bZ)^*$. So Galois action is
transitive on the non-zero classes in the kernel of $\hat\rho$.

Let $w_1$, $w_2$  and $w_3$ be three
non-colinear points in $X(\bK)$.
The linear pencil of the divisor $w_1+w_2+w_3$ has
projective dimension zero, that is  \[H^0(X,\cO_X(w_1+w_2+w_3))=\bK.\]
Let $v_1$, \dots,  $v_7$
be seven  points in $X(\bK)$.
 We assume that  $v_1$, $v_2$, $v_3$, $v_4$, $v_5$, $v_6$, $v_7$,
$w_1$, $w_2$ and $w_3$
are  pairwise distinct.
Since the cardinality of $J_X(\bK)$ is odd,
the multiplication by two map is a bijection of it.  We deduce the existence
 of seven    effective degree three divisors  $D_1$, \dots, $D_7$
such that
$2(D_i-3w_0)$ is linearly equivalent to $w_1+w_2+w_3-v_i-2w_0$ for
$1\leqslant i\leqslant 7$.
The divisor $2D_i-4w_0$ is linearly equivalent
to $w_1+w_2+w_3-v_i$. It is non-special.

Let $\xi$ be any non-zero
divisor class in the kernel of $\hat\rho$.
For each $1\leqslant i\leqslant 7$,
the divisor class $2D_i-4w_0-\xi$ is the class
$w_1+w_2+w_3-v_i-\xi$. At most  three
among these seven  classes are effective.
Otherwise the class $w_1+w_2+w_3-\xi$ would have positive projective dimension.
So it would be  the class  $K-P_\xi$ where
$K$ is the canonical class and $P_\xi$ is a point on $X$.
Because Galois action is transitive on the non-zero
classes in the kernel of $\hat\rho$,  there would exist for every
such class $\xi$ a point $P_\xi$ such that
$w_1+w_2+w_3-\xi$ is the class  $K-P_\xi$.
We consider the points $P_s$, $P_{2s}$, $P_{3s}$, $P_{4s}$
associated with  $s$, $2s$, $3s$, $4s$, where $s$ is a generator
of the kernel of $\hat\rho$. These are four pairwise distinct points
and $P_{2s}-P_s$ is linearly equivalent to $P_{4s}-P_{3s}$.
So the linear series of $P_{s}+P_{4s}$ has positive projective
dimension. A contradiction because $X$ is not hyperelliptic.

So we can assume that   $w_1+w_2+w_3-v_i-\xi$ is ineffective
for $1\leqslant i\leqslant 4$.
Since the Galois group of $\bK$ acts transitively on the non-zero
classes in the kernel of $\hat\rho$ we deduce that $2D_i -4w_0-\xi$
is ineffective for any $\xi$ in $\Ker\hat\rho$ and any $1\leqslant i\leqslant
4$.

At least one among $D_1$, $D_2$, $D_3$, $D_4$ has projective dimension zero.
Otherwise there would exist four  points $P_1$, $P_2$, $P_3$, $P_4$
such  that $D_i$ is linearly equivalent to $K-P_i$
for $1\leqslant i\leqslant 4$. So $2(K-P_i)$ is linearly
equivalent to $w_1+w_2+w_3+4w_0-v_i$. We deduce
that $2P_2+v_1 \sim 2P_1+v_2$ and $2P_3+v_1\sim
2P_1+v_3$ and $2P_4+v_1\sim
2P_1+v_4$. So these  classes have positive projective
dimension. There exist three  points $Q_{2}$, $Q_{3}$, and $Q_4$
such that $2P_1+v_2\sim K-Q_{2}$,  $2P_1+v_3\sim K-Q_{3}$,
and $2P_1+v_4\sim K-Q_{4}$. So $v_2+Q_{2}\sim v_3+Q_{3}\sim v_4+Q_4$.
A contradiction because $X$ is not hyperelliptic.

We call $D$  one among  $D_1$, $D_2$, $D_3$, $D_4$
having  projective dimension zero. And let $v$
be the point such that
$2(D-3w_0)$ is linearly equivalent to $w_1+w_2+w_3-v-2w_0$.

Because $n$ is a prime integer, every  fiber of $\rho$
above a rational point of $X$ is either irreducible
or completely split. Since de genus of $Y$ is
\[g_Y=1+n(g_X-1)=1+2n=9279,\] the number of $\bK$-rational
points on it is bounded from above by $q+1+2g_Y\sqrt q <66926$.
So the number of split fibers is $\leqslant 66925/9279<15$. So there are
at least $9$  points $(a_j)_{1\leqslant j\leqslant 9}$
in $X(\bK)$ with an irreducible
fiber above them.


Let $\xi$ be any non-zero class
in the kernel of $\hat\rho$.
At most three  among  the
$(a_i)_{1\leqslant i\leqslant 9}$
make $D-a_i-\xi$ effective. Otherwise $D-\xi$ would have
positive projective dimension. So it would be linearly
equivalent to $K-P_\xi$ for some
point $P_\xi$ on $X$.
Because Galois action is transitive on the non-zero
classes in the kernel of $\hat\rho$,  there would exist for every
such class $\xi$ a point $P_\xi$ such that
$D-\xi$ is the class  $K-P_\xi$.
We consider the points $P_s$, $P_{2s}$, $P_{3s}$, $P_{4s}$
associated with  $s$, $2s$, $3s$, $4s$, where $s$ is a generator
of the kernel of $\hat\rho$. These are four pairwise distinct points
and $P_{2s}-P_s$ is linearly equivalent to $P_{4s}-P_{3s}$.
So the linear series of $P_{s}+P_{4s}$ has positive projective
dimension. A contradiction because $X$ is not hyperelliptic.
So at least
six  among  the
$(a_i)_{1\leqslant i\leqslant 9}$
make $D-a_i-\xi$ ineffective for the chosen non-zero $\xi$
and thus for all its conjugates.
So we can assume
that 
$D-a_i-\xi$ is ineffective for any $1\leq i\leq 6$ and
any non-zero $\xi$ in the  kernel of $\hat\rho$.

At  most three  among  the
$(a_i)_{1\leqslant i\leqslant 6}$
make $D-a_i$ special. Otherwise $D$ would have
positive projective dimension. A contradiction.

We let $a$ be one among  $(a_i)_{1\leq i\leq 6}$
such that  $D-a$ is non-special.
We let  $\rgot$  be the divisor consisting
of four times any point on $Y$ above
$w_0$.
The ring $\bK[\rgot ]$ is isomorphic to $\bK[x]/x^4$. Since $\bK$ has $q\geqslant 6$
elements the symmetric bilinear complexity of
multiplication in $\bK[\rgot ]$ is $7$. Applying Theorem \ref{th:analy}
we deduce \[\nu^{sym}_{13}(4639)\leqslant 7.\]

\smallskip

The following  computer session implements
this calculation  in 
SageMath (Version 9.4)  \cite{sagemath}.

\smallskip

\begin{verbatim}
sage: q=13;K=GF(q);Kx.<x>=FunctionField(K);Kxy.<y>=Kx[];
KX.<y> = Kx.extension(y^3+y*(K(1)+K(2)*x+x^3)+K(10)*x
+K(9)*x^2+K(11)*x^3)
g=KX.genus();LP=KX.L_polynomial();t=LP.parent().gen();
sage: chi=LP(1/t)*t^(2*g)                                                                        
t^6 + 9*t^5 + 51*t^4 + 197*t^3 + 663*t^2 + 1521*t + 2197
sage: n=numerator(chi(1))
4639
Fnt.<t> = PolynomialRing(GF(n));factor(Fnt(chi))
(t+2195)*(t+3726)*(t+4626)*(t+4638)*(t^2+3380*t+13)
\end{verbatim}


\section{Remarks and questions}\label{sec:concl}

The symmetric equivariant complexity $\nu_q^{sym}(n)$
provides  a good control on  the computational
difficulty of multiplying two elements in a degree $n$ extension
of a field $\bK$ with $q$-elements, given by their coordinates in any
normal basis.
We have shown how to bound $\nu_q^{sym}(n)$ using points of order
$n$ in Jacobians over $\bK$. We  need  
a Jacobian with smallest possible  dimension having  a point of  order $n$.
A natural question is : given $q$ and $n$, which  is the smallest possible
$g$ such that there exists a Jacobian  of dimension $g$ over a field
with $q$ elements, having a rational point of order $n$ ?
Are there asymptotic families that are good
with this respect ?
Modular towers produce  curves with many points but they
have too much ramification to be useful here.

In practice, we pick random
curves of genus $g$ over a field with $q$ elements, until we find some
whose Jacobian has order divisible by $n$.
A difficulty is that for large $g$ we do not have
a convenient model for a universal  curve of genus $g$.
We could restrict to hyperelliptic curves but their
Jacobians tend to be smaller, so this restriction affects the
efficiency of the method.


 We may wonder if Theorem \ref{th:redu} is optimal, even roughly.
 Given $q$ and some bound  $C$, are there only finitely many $n$
 such that $\nu_q^{sym}(n)\leqslant C$ ?
 such that  $\nu_q^{sym}(n)\leqslant C|\log (\log (n))|$ for example ?

\bibliographystyle{plain}
\bibliography{revision}

\end{document}